\documentclass{amsart}
\usepackage{amsfonts}

 \theoremstyle{definition}
 
 \theoremstyle{remark}
 
 \numberwithin{equation}{subsection}

\begin{document}
\title{ ON CLASSIFICATION OF FINITE DIMENSIONAL COMPLEX FILIFORM LEIBNIZ
ALGEBRAS (Part 2)}
\author{Bekbaev U.D.}
\author{Rakhimov I.S.}
\address{Bekbaev Ural Djumaevich.\\
\indent
Department of Mathematics $\&$ Institute for Mathematical Research, FS,UPM,
43400, Serdang, Selangor Darul Ehsan, (Malaysia).}
\email{bekbaev@science.upm.edu.my}
\address{Rakhimov Isamiddin Sattarovich.\\
\indent
Department of Mathematics $\&$ Institute for Mathematical Research, FS,UPM,
43400, Serdang, Selangor Darul Ehsan, (Malaysia).}
\email{isamiddin@science.upm.edu.my \mbox{@} risamiddin@mail.ru}
\maketitle

\begin{abstract}
The paper is devoted to classification problem of finite
dimensional complex none Lie filiform Leibniz algebras. Actually,
the observations show there are two resources to get
classification of filiform Leibniz algebras. The first of them is
naturally graded none Lie filiform Leibniz algebras and the
another one is naturally graded filiform Lie algebras. Using the
first resource we get two disjoint classes of filiform Leibniz
algebras \cite{AyuB}. The present paper deals with the second of
the above two classes, the first class has been considered in
\cite{BR}. The algebraic classification here means to specify the
representatives of the orbits, whereas the geometric
classification is the problem of finding generic structural
constants in the sense of algebraic geometry. Our main effort in
this paper is the algebraic classification. We suggest here an
algebraic method based on invariants. Utilizing this method for
any given low dimensional case all filiform Leibniz algebras can
be classified. Moreover, the results can be used for geometric
classification of orbits of such algebras.
\end{abstract}



\section{Introduction}

This paper aims to investigate a class of nonassociative algebras
which generalizes the class of Lie algebras. These algebras
satisfy certain identities that were suggested by J.-L.Loday
\cite{Lod}. When he used the tensor product instead of external
product in the definition of the $n$-th cochain, in order to prove
the differential property, that is defined on cochains, it
sufficed to replace the anticommutativity and Jacoby identity by
the Leibniz identity. This is an essential one of the motivation
to appear for this class of algebras.

In this paper we suggest an algebraic approach to the
classification problem for filiform Leibniz algebras. Utilizing
this method for any fixed low dimensional case the corresponding
classes of filiform Leibniz algebras can be classified completely.
Moreover, the results may be used for geometric classification in
the sense of geometric invariant theory \cite{Kr}. It is assumed
that it will be the subject of one of the next papers. For
geometric classification of complex nilpotent Leibniz algebras of
dimension at most four we refer to \cite{AOR1}.

Let $V$ be a vector space of dimension $n$ over an algebraically closed
field $K$ (char$K$=0). The bilinear maps $V \times V \rightarrow V$ form a
vector space $Hom(V\otimes V,V)$ of dimension $n^{3}$, which can be
considered together with its natural structure of an affine algebraic
variety over $K$ and denoted by $Alg_{n}(K)\cong K^{n^{3}}$. An $n$%
-dimensional algebra $L$ over $K$ may be considered as an element $%
\lambda(L) $ of $Alg_{n}(K)$ via the bilinear mapping $\lambda:L\otimes L\to
L$ defining an binary algebraic operation on $L:$ let $\{
e_1,e_2,\ldots,e_n\}$ be a basis of the algebra $L.$ Then the table of
multiplication of $L$ is represented by point $(\gamma_{ij}^{k})$ of this
affine space as follow:
\begin{equation*}
\lambda(e_i,e_j)=\sum\limits_{k=1}^n\gamma_{ij}^ke_k.
\end{equation*}
$\gamma_{ij}^k$ are called \textit{structural constants} of $L.$ The linear
reductive group $GL_n(K)$ acts on $Alg_{n}(K)$ by $(g*\lambda)(x,y)=g(%
\lambda(g^{-1}(x),g^{-1}(y)))$("transport of structure"). Two algebras $%
\lambda_{1}$ and $\lambda_{2}$ are isomorphic if and only if they belong to
the same orbit under this action. The orbit of $\lambda$ under this action
is denoted by $O(\lambda)$. It is clear that elements of the given orbit are
isomorphic to each other algebras. The classification means to specify the
representatives of the orbits. A simple criterion, to decide if the given
two algebras are isomorphic, is desired.

\section{Preliminaries}

\textbf{Definition 1.} An algebra $L$ over a field $K$ is called a
\textit{Leibniz algebra} if it satisfies the following Leibniz
identity:
\begin{equation*}
[x,[y,z]]=[[x,y],z]-[[x,z],y],
\end{equation*}
where $[\cdot,\cdot]$ denotes the multiplication in $L.$ Let
$Leib_{n}(K)$ be a subvariety of $Alg_{n}(K)$ consisting of all
$n$-dimensional Leibniz algebras over K. It is invariant under the
above mentioned action of $GL_n(K)$. As a subset of $Alg_{n}(K)$
the set $Leib_{n}(K)$ is specified by system of equations with
respect to structural constants $\gamma_{ij}^{k}$:
\begin{equation*}
\sum\limits_{\emph{l}=1}^{\emph{n}}{(\gamma_{\emph{jk}}^{\emph{l}} \gamma_{%
\emph{il}}^{\emph{m}}-\gamma_{\emph{ij}}^{\emph{l}}\gamma_{\emph{lk}}^{\emph{%
m}}+ \gamma_{\emph{ik}}^{\emph{l}}\gamma_{\emph{lj}}^{\emph{m}})}=0
\end{equation*}
It is easy to see that if the bracket in Leibniz algebra happens
to be anticommutative then it is a Lie algebra. So Leibniz
algebras are "noncommutative" generalization of Lie algebras. As
to classifications of low dimensional Lie algebras they are well
known. But unless simple Lie algebras the classification problem
of all Lie algebras in common remains a big problem. Yu.I.Malcev
\cite{Mal} reduced the classification of solvable Lie algebras to
the classification of nilpotent Lie algebras. Apparently the first
non-trivial classification of some classes of low-dimensional
nilpotent Lie algebra are due to Umlauf. In his thesis \cite{Um}
he presented the redundant list of nilpotent Lie algebras of
dimension at most seven. He gave also the list of nilpotent Lie
algebras of dimension less than ten admitting so-called adapted
basis (now, the nilpotent Lie algebras with this property are
called \emph{filiform Lie algebras}). It was shown by M.Vergne
\cite{Vr} the importantness of filiform Lie algebras in the study
of variety of nilpotent Lie algebras laws. Up to now the several
classifications of low-dimensional nilpotent Lie algebras have
been offered. Unfortunately, many of these papers are based on
direct computations and the complexity of those computations leads
frequently to errors. We refer the reader to \cite{GJK} for
comments and corrections of the classification errors.

Further if it is not asserted additionally all algebras assumed to
be over the field of complex numbers $\textbf{C}.$

Let $L$ be a Leibniz algebra. We put:
\begin{equation*}
L^{1}=L,\quad L^{k+1}=[L^{k},L],k\in N.
\end{equation*}

\textbf{Definition 2.} A Leibniz algebra $L$ is said to be \textit{nilpotent}
if there exists an integer $s\in N,$ such that $L^{1} \supset L^{2} \supset
... \supset L^{s}=\{0\}.$ The smallest integer $s$ for which $L^{s}=0$ is
called \textit{the nilindex} of $L$.

\textbf{Definition 3.} An $n$-dimensional Leibniz algebra $L$ is
said to be \textit{filiform} if $dim L^i =n-i,$ for all $2\le i\le
n.$

\textbf{Theorem 1.}\cite{AyuB},\cite{GO}.

Any $(n+1)$-dimensional complex non Lie filiform Leibniz algebra
can be included to one of the following three classes of none Lie
filiform Leibniz algebras:
\newline
a) ($1^{\mbox{st}}$ \mbox{class}): \ \ $\left\{%
\begin{array}{ll}
[e_0,e_0]=e_{2}, &  \\[1mm]
[e_i,e_0]=e_{i+1}, & \ 1\leq i \leq {n-1} \\[1mm]
[e_0,e_1]=\alpha_3e_3 + \alpha_4e_4+...+ \alpha_{n-1}e_{n-1}+ \theta e_n, &
\\[1mm]
[e_j,e_1]=\alpha_3e_{j+2} + \alpha_4e_{j+3}+...+ \alpha_{n+1-j}e_n, & \
\ 1\leq j \leq {n-2}%
\end{array}
\right.$ \\[1mm]
(omitted products are supposed to be zero)

b) ($2^{\mbox{nd}}$ \mbox{class}): \ \ $\left\{%
\begin{array}{ll}
[e_0,e_0]=e_{2}, &  \\[1mm]
[e_i,e_0]=e_{i+1}, & \ 2\leq i \leq {n-1} \\[1mm]
[e_0,e_1]=\beta_3e_3 + \beta_4e_4+...+ \beta_ne_n, &  \\[1mm]
[e_1,e_1]=\gamma e_n, &  \\[1mm]
[e_j,e_1]=\beta_3e_{j+2} + \beta_4e_{j+3}+...+ \beta_{n+1-j}e_n, & \ 2\leq j
\leq {n-2}%
\end{array}
\right.$ \newline
(omitted products are supposed to be zero)

c) ($3^{\mbox{rd}}$ \mbox{class}):  $\left\{%
\begin{array}{lll}
[e_0,e_0]= e_{n}, &  &  \\[1mm]
[e_1,e_1]=\alpha e_{n}, &  &  \\[1mm]
[e_i,e_0]=e_{i+1}, & 1\leq i \leq {n-1} &  \\[1mm]
[e_0,e_1]=-e_2+\beta e_n, &  &  \\[1mm]
[e_0,e_i]=-e_{i+1}, & 2\leq i \leq {n-1} &  \\[1mm]
[e_i,e_j]= \\ \ \ -[e_j,e_i] \in lin<e_{i+j+1}, e_{i+j+2}, \dots ,
e_n>, & 1\leq i
\leq n-3, &  \\[1mm]
& 2 \leq j \leq {n-1-i} &  \\[1mm]
[e_{n-i},e_i]=-[e_i,e_{n-i}]=(-1)^i \delta e_n, & 1 \leq i \leq n-1 &  \\%
[1mm]
&  &
\end{array}
\right.$ \newline (omitted products are supposed to be zero)
\newline where $\{e_0, e_1, e_2, .... , e_n \}$ is a basis,
$\delta$ is either 1 or 0 for odd $n$ and $\delta=0$ for even $n.$

In other words, the above proposition means that the set of all
$(n+1)$-dimensional complex none Lie filiform Leibniz algebras can
be represented as a disjoint union of the above mentioned three
classes and the algebras from the difference classes never are
isomorphic to each other.

In this paper we will consider the second class of algebras.

Let us denote by $L(\beta)$, the $(n+1)$-dimensional filiform
non-Lie Leibniz algebra defined by parameters
$\beta=(\beta_{3},\beta_{4},...,\beta_{n},\gamma)$. The set of all
$(n+1)$-dimensional complex filiform Leibniz algebras from the
second class is denoted by $FLeib_{n+1}$. It is a closed and
invariant subset of the variety of nilpotent Leibniz algebras.

Using the method of simplification of the basis transformations in
\cite{GO} the following criterion on isomorphism of two
$(n+1)$-dimensional filiform Leibniz algebras was given. We
formulate part two of the theorem regarding our case. Namely: let
$n\geq 3$.

\textbf{Theorem 2.}\cite{GO} Two algebras $L(\beta_3, \beta_4, ...
, \beta_n, \gamma )$ and $L^{^{\prime}}(\beta^{\prime}_{3},
\beta^{\prime}_{4}, ... , \beta^{\prime}_{n}, \gamma^{\prime})$
from $FL_{n+1}$ are isomorphic if and only if there exist
$\emph{A},\emph{B} \ \ \mbox{and} \ \ \emph{D}\in \textbf{C}$ such
that $\emph{AD}\neq 0$ and the following conditions hold:
\begin{equation*}
{\small
\begin{array}{ll}

\beta^{\prime}_3=\frac{D}{A^2} \beta_3 &  \\
\beta^{\prime}_t=\frac{1}{A^{t-1}}(D\beta_t - \sum\limits_{k=3}^{t-1}(
C_{k-1}^{k-2}A^{k-2}B\beta_{t+2-k}+C_{k-1}^{k-3}A^{k-3}B^2\sum%
\limits_{i_1=k+2}^{t}\beta_{t+3-i_1} \cdot \beta_{i_1+1-k} + &  \\
\ \ \ \ +C_{k-1}^{k-4}A^{k-4}B^3\sum\limits_{i_2=k+3}^{t}
\sum\limits_{i_1=k+3}^{i_2}\beta_{t+3-i_2}\cdot
\beta_{i_2+3-i_1}\cdot
\beta_{i_1-k} + \ ... \ + &  \\
\ \ \ \ +C_{k-1}^{1}AB^{k-2}\sum\limits_{i_{k-3}=2k-2}^{t}
\sum\limits_{i_{k-4}=2k-2}^{i_{k-3}} ...
\sum\limits_{i_1=2k-2}^{i_2} \beta_{t+3-i_{k-3}}
\beta_{i_{k-3}+3-i_{k-4}} ...
\beta_{i_2+3-i_1}\beta_{i_1+5-2k}+ &  \\
\ \ \ \ + B^{k-1}\sum\limits_{i_{k-2}=2k-1}^{t}
\sum\limits_{i_{k-3}=2k-1}^{i_{k-2}} ....
\sum\limits_{i_1=2k-1}^{i_2} \beta_{t+3-i_{k-2}}
\beta_{i_{k-2}+3-i_{k-3}} ....
\beta_{i_2+3-i_1}\beta_{i_1+4-2k})\beta^{\prime}_k), &
\end{array}
}
\end{equation*}
where $4\leq t\leq n-1.$
\begin{equation*}
{\small
\begin{array}{ll}
\beta^{\prime}_n=\frac{BD\gamma}{A^n}+\frac{1}{A^{n-1}}(D\beta_n -
\sum\limits_{k=3}^{n-1}(
C_{k-1}^{k-2}A^{k-2}B\beta_{n+2-k}+C_{k-1}^{k-3}A^{k-3}B^2\sum%
\limits_{i_1=k+2}^{n}\beta_{n+3-i_1} \cdot \beta_{i_1+1-k} + &  \\
\ \ \ \ +C_{k-1}^{k-4}A^{k-4}B^3\sum\limits_{i_2=k+3}^{n}
\sum\limits_{i_1=k+3}^{i_2}\beta_{n+3-i_2}\cdot
\beta_{i_2+3-i_1}\cdot
\beta_{i_1-k} + \ ... \ + &  \\
\ \ \ \ +C_{k-1}^{1}AB^{k-2}\sum\limits_{i_{k-3}=2k-2}^{n}
\sum\limits_{i_{k-4}=2k-2}^{i_{k-3}} ...
\sum\limits_{i_1=2k-2}^{i_2} \beta_{n+3-i_{k-3}}
\beta_{i_{k-3}+3-i_{k-4}} ...
\beta_{i_2+3-i_1}\beta_{i_1+5-2k}+ &  \\
\ \ \ \ + B^{k-1}\sum\limits_{i_{k-2}=2k-1}^{n}
\sum\limits_{i_{k-3}=2k-1}^{i_{k-2}} ....
\sum\limits_{i_1=2k-1}^{i_2} \beta_{n+3-i_{k-2}}
\beta_{i_{k-2}+3-i_{k-3}} ....
\beta_{i_2+3-i_1}\beta_{i_1+4-2k})\beta^{\prime}_k), \\

\gamma^{\prime}=\frac{D^2}{A^{n}} \gamma &  \\[1mm]
\end{array}
}
\end{equation*}

Here are the above systems of equalities for some low dimensional cases:

Case of $n=4$ i.e. $\mathbf{dimL=5:}$
\begin{equation}
\begin{array}{lll}
\left\{
\begin{array}{lll}
\beta^{\prime}_{3}=\frac{1}{A}\frac{D}{A}\beta_{3},\\ \\
\beta^{\prime}_{4}=\frac{1}{A^2}\frac{D}{A}(\frac{B}{A}\gamma+\beta_{4}-
2\frac{B}{A}\beta_{3}^{2}).\\ \\
\gamma^{\prime}=\frac{1}{A^2}(\frac{D}{A})^{2}\gamma, \\

\end{array}
\right.
\end{array}
\end{equation}

Case of $n=5$ i.e. $\mathbf{dimL=6:}$
\begin{equation}
\begin{array}{lll}
\left\{
\begin{array}{lll}
\beta^{\prime}_{3}=\frac{1}{A}\frac{D}{A}\beta_{3},\\ \\
\ \beta^{\prime}_{4}=\frac{1}{A^2}\frac{D}{A}
(\beta_{4}-2\frac{B}{A}\beta_{3}^{2}), \\
\\
\beta^{\prime}_{5}=\frac{1}{A^{3}}\frac{D}{A}[\frac{B}{A}\gamma+
\beta_{5}-5\frac{B}{A}\beta_3\beta_{4}+5(\frac{B}{A})^2\alpha_{3}^{3}].
\\ \\
\gamma^{\prime}=\frac{1}{A^3}(\frac{D}{A})^{2}\gamma, \\
\end{array}
\right.
\end{array}
\end{equation}

Case of $n=6$ i.e. $\mathbf{dimL=7:}$
\begin{equation}
\begin{array}{lll}
\left\{
\begin{array}{lll}
\beta^{\prime}_{3}=\frac{1}{A}\frac{D}{A}\beta_{3}, \\ \\

\beta^{\prime}_{4}=\frac{1}{A^2}\frac{D}{A}(\beta_{4}-2\frac{B}{A}\beta_{3}^{2}), \\
\\
\beta^{\prime}_{5}=\frac{1}{A^{3}}\frac{D}{A}[(\beta_{5}-5\frac{B}{A}
\beta_3\beta_{4}+5(\frac{B}{A})^2\beta_{3}^{3}]. \\
\\
\beta^{\prime}_{6}=\frac{1}{A^{4}}\frac{D}{A}[\frac{B}{A}\gamma+\beta_6-6
\frac{B}{A}\beta_3\beta_5+
21(\frac{B}{A})^2\beta_3^2\beta_4-3\frac{B}{A}
\beta_4^2-14(\frac{B}{A})^3\beta_{3}^{4}]. \\
\\
\gamma^{\prime}=\frac{1}{A^4}(\frac{D}{A})^{2}\gamma, \\ \\
\end{array}
\right.
\end{array}
\end{equation}

Case of $n=7$ i.e. $\mathbf{dimL=8:}$
\begin{equation}
\begin{array}{lll}
\left\{
\begin{array}{lll}
\beta^{\prime}_{3}=\frac{1}{A}\frac{D}{A}\beta_{3}, \\ \\

\beta^{\prime}_{4}=\frac{1}{A^2}\frac{D}{A}(\beta_{4}-2\frac{B}{A}\beta_{3}^{2}), \\
 \\
\beta^{\prime}_{5}=\frac{1}{A^{3}}\frac{D}{A}[(\beta_{5}-5\frac{B}{A}
\beta_3\beta_{4}+5(\frac{B}{A})^2\beta_{3}^{3}]. \\
\\
\beta^{\prime}_{6}=\frac{1}{A^{4}}\frac{D}{A}[\beta_6-6\frac{B}{A}
\beta_3\beta_5+
21(\frac{B}{A})^2\beta_3^2\beta_4-3\frac{B}{A}\beta_4^2-14(
\frac{B}{A})^3\beta_{3}^{4}] \\
 \\
\beta^{\prime}_{7}=\frac{1}{A^{5}}\frac{D}{A}[\frac{B}{A}\gamma+\beta_{7}-7
\frac{B}{A}\beta_3\beta_6+
28(\frac{B}{A})^2\beta_3^2\beta_5+28(\frac{B}{A}
)^2\beta_3\beta_4^2- \\ \\ \qquad \qquad \qquad
\qquad \qquad \qquad \qquad 7\frac{B}{A}\beta_4\beta_5-84(\frac{B}{A})^3\beta_{3}^{3}%
\beta_4+42(\frac{B}{A})^4\beta_{3}^{5}]
\\ \\
\gamma^{\prime}=\frac{1}{A^5}(\frac{D}{A})^{2}\gamma, \\ \\
\end{array}
\right.
\end{array}
\end{equation}

To deal with the classification of $FLeib_{n+1}$ with respect to
the above mentioned action we represent it as a disjoint union of
an open and closed (with respect to the Zarisski topology)
subsets. Moreover each of these subsets are invariant under the
corresponding transformations presented in Theorem 2. Then we
formulate the solution of the isomorphism problem for the
corresponding algebras from the open subset. Similar approach can
be used to solve isomorphism problem for the algebras from the
corresponding closed subset.

It is not difficult to notice that the expressions for $\beta^{\prime}_{t}$,
$\gamma^{\prime}$ in Theorem 2 can be represented in the following form:
\begin{equation}
\beta^{\prime}_{t}=\frac{1}{A^{t-2}}\frac{D}{A}\psi_{t}(\frac{B}{A};\beta),
\end{equation}
where $\beta=(\beta_{3},\beta_{4},..., \beta_{\emph{n}},\gamma)$, $3 \leq t
\leq n-1$ and

$\psi_{t}(y;z)=\psi_{t}(y;z_{3},z_{4},...,z_{n},z_{n+1})= z_{t}- \sum
\limits_{k=3}^{t-1}(C_{k-1}^{k-2}yz_{t+2-k}+ C_{k-1}^{k-3}y^{2} \sum
\limits_{i_{1}=k+2}^{t} z_{t+3-i_{1}}\cdot
z_{i_{1}+1-k}+C_{k-1}^{k-4}y^{3}\sum\limits_{i_{2}=k+3}^{t}
\sum\limits_{i_{1}=k+3}^{i_{2}}z_{t+3-i_{2}} \cdot z_{i_{2}+3-i_{1}} \cdot
z_{i_{1}-k}+...+$\newline
$C_{k-1}^{1}y^{k-2}\sum\limits_{i_{k-3}=2k-2}^{t}
\sum\limits_{i_{k-4}=2k-2}^{i_{k-3}}...
\sum\limits_{i_{1}=2k-2}^{i_{2}}z_{t+3-i_{k-3}} \cdot z_{i_{k-3}+3-i_{k-4}}
\cdot ...\cdot z_{i_{2}+3-i_{1}} \cdot z_{i_{1}+5-2k}$\newline
$+y^{k-1} \sum\limits_{i_{k-2}=2k-1}^{t}
\sum\limits_{i_{k-3}=2k-1}^{i_{k-2}}...\sum%
\limits_{i_{1}=2k-1}^{i_{2}}z_{t+3-i_{k-2}} \cdot
z_{i_{k-2}+3-i_{k-3}} \cdot...\cdot z_{i_{2}+3-i_{1}} \cdot
z_{i_{1}+4-2k})\cdot \psi_{k}(y;z),$

for $3 \leq t \leq n,$
\begin{equation}
\begin{array}{ll}
\beta^{\prime}_{n}=\frac{1}{A^{n-2}}\frac{D}{A}\frac{B}{A}\gamma+\psi_{n}(
\frac{B}{A};\beta),
\end{array}
\end{equation}
and
\qquad
$\gamma^{\prime}=\frac{1}{A^{n-2}}(\frac{D}{A})^2\psi_{n+1}(\frac{A}{B
};\beta),$ where $\psi_{n+1}(y;z)=z_{n+1}$

To simplify notation let us agree that in the above case for
transition from the $(n+1)$-dimensional filiform Leibniz algebra
$L(\beta)$ to the $n+1$-dimensional filiform Leibniz algebra
$L(\beta^{\prime})$\ \ we write
$\beta^{\prime}=\varrho(\frac{1}{A},\frac{B}{A},\frac{D}{A};\beta),$
\ \ where $\beta=(\beta_{3},\beta_{4},...,\beta_{n},\gamma)$

\begin{center}
$\varrho(\frac{1}{A},\frac{B}{A},\frac{D}{A};\beta)=
(\varrho_{1}(\frac{1}{A},\frac{B}{A},\frac{D}{A};\beta),
\varrho_{2}(\frac{1}{A},\frac{B}{A},\frac{D}{A};\beta),
...,\varrho_{n-1}(\frac{1}{A},\frac{B}{A},\frac{D}{A};\beta)),$ \
\ $\varrho_{t}(x,y,u;z)=x^{t-1}u\psi_{t+2}(y;z)$ \ \ for \ \ $1
\leq t \leq n-2$ and $\varrho_{n-1}(x,y,u;\overline{z})=
x^{n-5}u^{2}\psi_{n+1}(y;z)$
\end{center}

Here are the main properties, used in this paper, of the operator
$\varrho$:
\newline

$1^{0}. \ \ \varrho(1,0,1;\cdot)) \ \ \mbox{is the
identity operator}.$\newline

$2^{0}. \ \
\varrho(\frac{1}{A_{2}},\frac{B_{2}}{A_{2}},\frac{D_{2}}{A_{2}};
\varrho(\frac{1}{A_{1}},\frac{B_{1}}{A_{1}},
\frac{D_{1}}{A_{1}};\beta))=
\varrho(\frac{1}{A_{1}A_{2}},\frac{B_{1}A_{2}+B_{2}D_{1}}{A_{1}A_{2}},
\frac{D_{1}D_{2}}{A_{1}A_{2}};\beta) $\newline

$3^{0} \ \ \mbox{If} \ \
\beta^{\prime}=\varrho(\frac{1}{A},\frac{B}{A}, \frac{D}{A};\beta)
\ \ \mbox{then} \ \ \beta=\varrho(A,-\frac{B}{D},
\frac{A}{D};\beta^{\prime})$.

\section{Classification theorems}

\textbf{Definition 4.} An action of algebraic group $G$ on a variety $Z$ is
a morphism

$\sigma: G \times Z \longrightarrow Z$ with

$(i)\ \ \sigma(e,z)=z$, where $e$ is the unit element of $G$ and
$z \in Z.$

$(ii)\ \ \sigma(g,\sigma(h,z))=\sigma(gh,z)$, for any $g,h \in G$
and $z \in Z.$

We shortly write $gz$ for $\sigma(g,z)$, and call $Z$ a
$G$-variety.

\textbf{Definition 5.} A morphism $\emph{f}:Z \longrightarrow K$,
($K$ is a base field) is said to be invariant if $\emph{f}(gz)=
\emph{f}(z)$ for any $g \in G$ and $z \in Z.$

The algebra of invariant morphisms on $Z$ with respect to the
action of the group G is denoted by $K[Z]^{G}$. Sometimes this
algebra is a finitely generated $K$-algebra. This is referred to
in \cite{W} as the "first fundamental problem of invariant
theory". If $Z$ is an irreducible then the field of rational
invariants can be defined as a quotient field of $K[Z]^{G}$. It is
always finitely generated as a subalgebra of the finitely
generated algebra $K(Z)$. Description the field of rational
invariants is an another important classical problem of the
invariant theory \cite{VP}.

Actually, we use some elements of the algebra of invariant morphisms under
the above mentioned adapted action on the variety of filiform Leibniz
algebras to solve isomorphism problem.

From here on we assume that $n\geq 5$ is a positive integer, since
there are complete classifications of complex nilpotent Leibniz
algebras of dimension at most four \cite{O},\cite{AOR2} (for
five-dimensional case see section 4).

Consider the following presentation of $FLeib_{n+1}$: \newline
\begin{equation}
FLeib_{n+1}=\emph{U}\cup \emph{F},
\end{equation}
\newline
where $\emph{U}= \{L(\beta):
\beta_{3}(4\beta_{3}^2\beta_6-12\beta_3\beta_4\beta_6+\beta_{4}^3)(4\beta_{3}
\beta_5-5\beta_4^2) \neq 0\}$, $\emph{F}= \{L(\beta):
\beta_{3}(4\beta_{3}^2\beta_6-12\beta_3\beta_4\beta_6+\beta_{4}^3)(4\beta_{3}
\beta_5-5\beta_4^2) = 0\}.$

Our main interest will be the cases of open sets, "generic
algebras", cases.

\textbf{Theorem 3.} $i)$ Two algebras $L(\beta)$ and $L(\beta')$
from $\emph{U}$ are isomorphic if and only if

$
\varrho_{\emph{i}}(\frac{\beta_{3}(4\beta_{3}\beta_{5}-5\beta_{4}^2)}{4(\beta_{3}^{2}\beta_{6}-
3\beta_{3}\beta_{4}\beta_{5}+2\beta_{4}^{3})}, \frac{\beta_{4}}{2
\beta_{3}^{2}},\frac{4(\beta_{3}^{2}\beta_{6}-3\beta_{3}\beta_{4}\beta_{5}+2\beta_{4}^{3})}
{\beta_{3}^2(4\beta_{3}\beta_{5}-5\beta_{4}^2)};\beta)$ \\ \\

$\qquad \qquad \qquad \qquad =
\varrho_{\emph{i}}(\frac{\beta_{3}'(4\beta_{3}'\beta_{5}'-5\beta_{4}'^2)}{4(\beta_{3}'^{2}\beta_{6}'-
3\beta_{3}'\beta_{4}'\beta_{5}'+2\beta_{4}'^{3})},\frac{\beta_{4}'}{2\beta_{3}'^{2}},\frac{4(\beta_{3}'^{2}\beta_{6}'-
3\beta_{3}'\beta_{4}'\beta_{5}'+2\beta_{4}'^{3})}{\beta_{3}'^2(4\beta_{3}'\beta_{5}'-5\beta_{4}'^2)};\beta')$
\\

whenever \ \ $i=\overline{3,n-1}.$

$ii)$ For any $(\lambda_3,\lambda_4,...,\lambda_{n-1})\in
\textbf{C}^{n-3}$ there is an algebra $L(\beta)$ from $\emph{U}$
such that
\begin{equation*}
\varrho_{\emph{i}}(\frac{\beta_{3}(4\beta_{3}\beta_{5}-5\beta_{4}^2)}{4(\beta_{3}^{2}\beta_{6}-
3\beta_{3}\beta_{4}\beta_{5}+2\beta_{4}^{3})}, \frac{\beta_{4}}{2
\beta_{3}^{2}},\frac{4(\beta_{3}^{2}\beta_{6}-3\beta_{3}\beta_{4}\beta_{5}+2\beta_{4}^{3})}
{\beta_{3}^2(4\beta_{3}\beta_{5}-5\beta_{4}^2)};\beta)= \lambda_i
\ \ \mbox{for all} \ \ \emph{i}=\overline{3,n-1}
\end{equation*}

\textbf{Proof.} $i).$ Let first two algebras $L(\beta )$ and
$L(\beta')$ be isomorphic that is to say there exist
$\emph{A,B,D}\in \textbf{C}$ such that $\emph{AD}\neq 0$ and
$\beta'=\varrho (\frac{1}{A},\frac{B}{A},\frac{D}{A};\beta)$.
Consider algebra $L(\beta^{0}),$ where $\beta ^{0}=\varrho
(\frac{1}{A_{0}},\frac{B_{0}}{A_{0}},\frac{D_{0}}{A_{0}};\beta )$
and $ A_{0}=\frac{4(\beta_{3}^{2}\beta_{6}-
3\beta_{3}\beta_{4}\beta_{5}+2\beta_{4}^{3})}{\beta_{3}(4\beta_{3}\beta_{5}-5\beta_{4}^2)},$
\newline $B_{0}=\frac{2\beta_{4}(\beta_{3}^{2}\beta_{6}-
3\beta_{3}\beta_{4}\beta_{5}+2\beta_{4}^{3})}{\beta_{3}^3(4\beta_{3}\beta_{5}-5\beta_{4}^2)}$
and $D_{0}=\frac{4(\beta_{3}^{2}\beta_{6}-
3\beta_{3}\beta_{4}\beta_{5}+2\beta_{4}^{3})}{\beta_{3}^2(4\beta_{3}\beta_{5}-5\beta_{4}^2)}.$
Since $\beta =\varrho (A,\frac{-B}{D},\frac{A}{D};\beta')$ and
$\beta ^{0}=\varrho
(\frac{1}{A_{0}},\frac{B_{0}}{A_{0}},\frac{D_{0}}{A_{0}};\beta
)=\varrho (\frac{1}{A_{0}},\frac{B_{0}}{A_{0}},
\frac{D_{0}}{A_{0}};\varrho
(A,\frac{-B}{D},\frac{A}{D};\beta'))\newline =\varrho
(\frac{A}{A_{0}},\frac{B_{0}A-A_{0}B}{A_{0}D},\frac{D_{0}A}{A_{0}D};\beta').$
It is easy to check that
$\frac{A}{A_{0}}=\frac{\beta_{3}'(4\beta_{3}'\beta_{5}'-5\beta_{4}'^2)}{4(\beta_{3}'^{2}\beta_{6}'-
3\beta_{3}'\beta_{4}'\beta_{5}'+2\beta_{4}'^{3})},$
$\frac{B_{0}A-A_{0}B}{A_{0}D}=\frac{\beta _{4}^{\prime }}{2\beta
_{3}'^{2}},$ and
$\frac{D_{0}A}{A_{0}D}=\frac{4(\beta_{3}'^{2}\beta_{6}'-
3\beta_{3}'\beta_{4}'\beta_{5}'+2\beta_{4}'^{3})}{\beta_{3}'^2(4\beta_{3}'\beta_{5}'-5\beta_{4}'^2)}.$

Therefore

$
\varrho(\frac{\beta_{3}(4\beta_{3}\beta_{5}-5\beta_{4}^2)}{4(\beta_{3}^{2}\beta_{6}-
3\beta_{3}\beta_{4}\beta_{5}+2\beta_{4}^{3})}, \frac{\beta_{4}}{2
\beta_{3}^{2}},\frac{4(\beta_{3}^{2}\beta_{6}-3\beta_{3}\beta_{4}\beta_{5}+2\beta_{4}^{3})}
{\beta_{3}^2(4\beta_{3}\beta_{5}-5\beta_{4}^2)};\beta)$ \\ \\

$\qquad \qquad \qquad \qquad =
\varrho(\frac{\beta_{3}'(4\beta_{3}'\beta_{5}'-5\beta_{4}'^2)}{4(\beta_{3}'^{2}\beta_{6}'-
3\beta_{3}'\beta_{4}'\beta_{5}'+2\beta_{4}'^{3})},\frac{\beta_{4}'}{2\beta_{3}'^{2}},\frac{4(\beta_{3}'^{2}\beta_{6}'-
3\beta_{3}'\beta_{4}'\beta_{5}'+2\beta_{4}'^{3})}{\beta_{3}'^2(4\beta_{3}'\beta_{5}'-5\beta_{4}'^2)};\beta')$

and, in particular,

$ \varrho
_{\emph{i}}(\frac{\beta_{3}(4\beta_{3}\beta_{5}-5\beta_{4}^2)}{4(\beta_{3}^{2}\beta_{6}-
3\beta_{3}\beta_{4}\beta_{5}+2\beta_{4}^{3})}, \frac{\beta_{4}}{2
\beta_{3}^{2}},\frac{4(\beta_{3}^{2}\beta_{6}-3\beta_{3}\beta_{4}\beta_{5}+2\beta_{4}^{3})}
{\beta_{3}^2(4\beta_{3}\beta_{5}-5\beta_{4}^2)};\beta)$ \\ \\

$\qquad \qquad \qquad \qquad = \varrho
_{\emph{i}}(\frac{\beta_{3}'(4\beta_{3}'\beta_{5}'-5\beta_{4}'^2)}{4(\beta_{3}'^{2}\beta_{6}'-
3\beta_{3}'\beta_{4}'\beta_{5}'+2\beta_{4}'^{3})},\frac{\beta_{4}'}{2\beta_{3}'^{2}},\frac{4(\beta_{3}'^{2}\beta_{6}'-
3\beta_{3}'\beta_{4}'\beta_{5}'+2\beta_{4}'^{3})}{\beta_{3}'^2(4\beta_{3}'\beta_{5}'-5\beta_{4}'^2)};\beta')$

for all $i=\overline{3,n-1}.$

This procedure can be shown schematically by the following
picture:

$\qquad \qquad  \qquad \qquad \qquad \beta \qquad
\buildrel(\frac{1}{A_{0}},\frac{B_{0}}{A_{0}},\frac{D_{0}}{A_{0}})\over
\longrightarrow  \ \quad  \beta^{0}$

$\qquad \qquad \qquad \qquad
(\frac{1}{A},\frac{B}{A},\frac{D}{A})\searrow \qquad \qquad
\nearrow
(\frac{A}{A_{0}},\frac{B_{0}A-A_{0}B}{A_{0}D},\frac{D_{0}A}{A_{0}D})$

\qquad \qquad

$\qquad \qquad \qquad \qquad \qquad \qquad \qquad \qquad \beta^{\prime}$

\bigskip

\bigskip Conversely, let the equalities  \\ \newline

$ \varrho
_{\emph{i}}(\frac{\beta_{3}(4\beta_{3}\beta_{5}-5\beta_{4}^2)}{4(\beta_{3}^{2}\beta_{6}-
3\beta_{3}\beta_{4}\beta_{5}+2\beta_{4}^{3})}, \frac{\beta_{4}}{2
\beta_{3}^{2}},\frac{4(\beta_{3}^{2}\beta_{6}-3\beta_{3}\beta_{4}\beta_{5}+2\beta_{4}^{3})}
{\beta_{3}^2(4\beta_{3}\beta_{5}-5\beta_{4}^2)};\beta)$ \\ \\

$\qquad \qquad \qquad \qquad = \varrho
_{\emph{i}}(\frac{\beta_{3}'(4\beta_{3}'\beta_{5}'-5\beta_{4}'^2)}{4(\beta_{3}'^{2}\beta_{6}'-
3\beta_{3}'\beta_{4}'\beta_{5}'+2\beta_{4}'^{3})},\frac{\beta_{4}'}{2\beta_{3}'^{2}},\frac{4(\beta_{3}'^{2}\beta_{6}'-
3\beta_{3}'\beta_{4}'\beta_{5}'+2\beta_{4}'^{3})}{\beta_{3}'^2(4\beta_{3}'\beta_{5}'-5\beta_{4}'^2)};\beta')$
\\

hold for \ \ $i=\overline{3,n-1}.$ Then it is easy to see that \\

$ \varrho
_{\emph{i}}(\frac{\beta_{3}(4\beta_{3}\beta_{5}-5\beta_{4}^2)}{4(\beta_{3}^{2}\beta_{6}-
3\beta_{3}\beta_{4}\beta_{5}+2\beta_{4}^{3})}, \frac{\beta_{4}}{2
\beta_{3}^{2}},\frac{4(\beta_{3}^{2}\beta_{6}-3\beta_{3}\beta_{4}\beta_{5}+2\beta_{4}^{3})}
{\beta_{3}^2(4\beta_{3}\beta_{5}-5\beta_{4}^2)};\beta)$ \\ \\

$\qquad \qquad \qquad \qquad = \varrho
_{\emph{i}}(\frac{\beta_{3}'(4\beta_{3}'\beta_{5}'-5\beta_{4}'^2)}{4(\beta_{3}'^{2}\beta_{6}'-
3\beta_{3}'\beta_{4}'\beta_{5}'+2\beta_{4}'^{3})},\frac{\beta_{4}'}{2\beta_{3}'^{2}},\frac{4(\beta_{3}'^{2}\beta_{6}'-
3\beta_{3}'\beta_{4}'\beta_{5}'+2\beta_{4}'^{3})}{\beta_{3}'^2(4\beta_{3}'\beta_{5}'-5\beta_{4}'^2)};\beta')$
\\

for $\ \emph{i}=\overline{1,2}$ as well and therefore

$$
\varrho(\frac{\beta_{3}(4\beta_{3}\beta_{5}-5\beta_{4}^2)}{4(\beta_{3}^{2}\beta_{6}-
3\beta_{3}\beta_{4}\beta_{5}+2\beta_{4}^{3})}, \frac{\beta_{4}}{2
\beta_{3}^{2}},\frac{4(\beta_{3}^{2}\beta_{6}-3\beta_{3}\beta_{4}\beta_{5}+2\beta_{4}^{3})}
{\beta_{3}^2(4\beta_{3}\beta_{5}-5\beta_{4}^2)};\beta)$$ \\ \\
$$\qquad \qquad \qquad \qquad =
\varrho(\frac{\beta_{3}'(4\beta_{3}'\beta_{5}'-5\beta_{4}'^2)}{4(\beta_{3}'^{2}\beta_{6}'-
3\beta_{3}'\beta_{4}'\beta_{5}'+2\beta_{4}'^{3})},\frac{\beta_{4}'}{2\beta_{3}'^{2}},\frac{4(\beta_{3}'^{2}\beta_{6}'-
3\beta_{3}'\beta_{4}'\beta_{5}'+2\beta_{4}'^{3})}{\beta_{3}'^2(4\beta_{3}'\beta_{5}'-5\beta_{4}'^2)};\beta')$$

which means the algebras $L(\beta )$ and $L(\beta')$ are
isomorphic to the same algebra and therefore they are isomorphic
to each other.

Part $ii)$ can be proved in the same way as the proof of Theorem 3
of \cite{BR}.

Here are the corresponding invariants for low dimensional cases.

Case of \textbf{dim L=6:}

$\ \ \ \ \ \ \ \ \ \ \ \ \ \ \ \ \ \ \ \ \ \ \ \ \ \ \ \ \ \ \ \ \
\ \ \varrho _{3}(\frac{\beta _{3}^{2}}{\gamma },\frac{\beta
_{4}}{2\beta _{3}^{2}},\frac{1}{\beta _{3}^{3}};\beta)=\frac{\beta
_{3}(4\beta _{5}\beta _{3}^{2}-5\beta _{4}^{2}\beta _{3}+2\beta
_{4}\gamma )}{4\gamma ^{2}}$

Case of \textbf{dim L=7:}\\ \\
\bigskip $
\begin{array}{ll}
\varrho _{3}(\frac{\beta
_{3}(4\beta_{3}\beta_{5}-5\beta_{4}^2)}{2(2\beta_{3}^2\beta_{6}-6\beta_{3}\beta_{4}\beta_{5}+\beta_{4}\gamma+4\beta_{4}^3)},
\frac{\beta _{4}}{2\beta
_{3}^{2}},\frac{2(2\beta_{3}^2\beta_{6}-6\beta_{3}\beta_{4}\beta_{5}+\beta_{4}\gamma+4\beta_{4}^3)}{\beta
_{3}^2(4\beta_{3}\beta_{5}-5\beta_{4}^2)};\beta)= \\ \\
\varrho_{4}(\frac{\beta
_{3}(4\beta_{3}\beta_{5}-5\beta_{4}^2)}{2(2\beta_{3}^2\beta_{6}-6\beta_{3}\beta_{4}\beta_{5}+\beta_{4}\gamma+4\beta_{4}^3)},
\frac{\beta _{4}}{2\beta
_{3}^{2}},\frac{2(2\beta_{3}^2\beta_{6}-6\beta_{3}\beta_{4}\beta_{5}+\beta_{4}\gamma+4\beta_{4}^3)}{\beta
_{3}^2(4\beta_{3}\beta_{5}-5\beta_{4}^2)};\beta) \\ \qquad \qquad
\qquad \qquad \qquad \qquad \qquad \qquad \qquad \qquad =
\frac{(4\beta _{3}\beta _{5}-5\beta _{4}^{2})^{3}}{16(\beta
_{4}\gamma +2\beta _{3}^{2}\beta _{6}-6\beta _{3}\beta _{4}\beta
_{5}+4\beta _{4}^{3})^{2}}
\end{array}
$

$\varrho _{5}(\frac{\beta _{3}^{2}}{\gamma },\frac{\beta
_{4}}{2\beta _{3}^{2}},\frac{1}{\beta _{3}^{3}};\beta
)=\frac{\gamma(4\beta _{3}\beta _{5}-5\beta _{4}^{2})^{2}}{4(\beta
_{4}\gamma +2\beta _{3}^{2}\beta _{6}-6\beta _{3}\beta _{4}\beta
_{5}+4\beta _{4}^{3})^{2}}$
\\

Case of \textbf{dim L=8:}\\ \\
\bigskip $
\begin{array}{ll}
\varrho _{3}(\frac{\beta
_{3}(4\beta_{3}\beta_{5}-5\beta_{4}^2)}{2(2\beta_{3}^2\beta_{6}-6\beta_{3}\beta_{4}\beta_{5}+\beta_{4}\gamma+4\beta_{4}^3)},
\frac{\beta _{4}}{2\beta
_{3}^{2}},\frac{2(2\beta_{3}^2\beta_{6}-6\beta_{3}\beta_{4}\beta_{5}+\beta_{4}\gamma+4\beta_{4}^3)}{\beta
_{3}^2(4\beta_{3}\beta_{5}-5\beta_{4}^2)};\beta) \\ \\ =
\varrho_{4}(\frac{\beta
_{3}(4\beta_{3}\beta_{5}-5\beta_{4}^2)}{2(2\beta_{3}^2\beta_{6}-6\beta_{3}\beta_{4}\beta_{5}+\beta_{4}\gamma+4\beta_{4}^3)},
\frac{\beta _{4}}{2\beta
_{3}^{2}},\frac{2(2\beta_{3}^2\beta_{6}-6\beta_{3}\beta_{4}\beta_{5}+\beta_{4}\gamma+4\beta_{4}^3)}{\beta
_{3}^2(4\beta_{3}\beta_{5}-5\beta_{4}^2)};\beta) \\ \qquad \qquad
\qquad \qquad \qquad \qquad \qquad \qquad \qquad \qquad =
\frac{(4\beta _{3}\beta _{5}-5\beta _{4}^{2})^{3}}{16(\beta
_{4}\gamma +2\beta _{3}^{2}\beta _{6}-6\beta _{3}\beta _{4}\beta
_{5}+4\beta _{4}^{3})^{2}}
\end{array}
$

\bigskip $\varrho _{5}(\frac{\beta
_{3}(4\beta_{3}\beta_{5}-5\beta_{4}^2)}{2(2\beta_{3}^2\beta_{6}-6\beta_{3}\beta_{4}\beta_{5}+\beta_{4}\gamma+4\beta_{4}^3)},
\frac{\beta _{4}}{2\beta
_{3}^{2}},\frac{2(2\beta_{3}^2\beta_{6}-6\beta_{3}\beta_{4}\beta_{5}+\beta_{4}\gamma+4\beta_{4}^3)}{\beta
_{3}^2(4\beta_{3}\beta_{5}-5\beta_{4}^2)};\beta)$ \\

$\qquad \qquad \qquad \qquad \qquad \qquad \qquad = \frac{(4\beta
_{3}\beta _{5}-5\beta
_{4}^{2})^{4}(4\beta_{3}\beta_{4}\gamma+8\beta_{3}^3\beta_{7}-28\beta_{3}^2\beta_{4}\beta_{6}
+28\beta_{3}\beta_{4}^2\beta_{5}-7\beta_{4}^4)}{128(2\beta_{3}^2\beta_{6}-6\beta_{3}\beta
_{4}\beta _{5}+4\beta _{4}^{3})}$

\bigskip $
\begin{array}{lll}
\varrho _{6}(\frac{\beta
_{3}(4\beta_{3}\beta_{5}-5\beta_{4}^2)}{2(2\beta_{3}^2\beta_{6}-6\beta_{3}\beta_{4}\beta_{5}+\beta_{4}\gamma+4\beta_{4}^3)},
\frac{\beta _{4}}{2\beta
_{3}^{2}},\frac{2(2\beta_{3}^2\beta_{6}-6\beta_{3}\beta_{4}\beta_{5}+\beta_{4}\gamma+4\beta_{4}^3)}{\beta
_{3}^2(4\beta_{3}\beta_{5}-5\beta_{4}^2)};\beta)=\frac{\beta
_{3}\gamma(4\beta _{3}\beta _{5}-5\beta _{4}^{2})^{3}}{8(2\beta
_{3}^{2}\beta _{6}-6\beta _{3}\beta _{4}\beta _{5}+4\beta
_{4}^{3})^{3}}
\end{array}
$
\\ \\

As to the isomorphism problem for the algebras from the closed set
$\emph{F}$ a similar procedure as the above can be applied to it
as well. We will not consider it here. In the next section we will
present final results of classification in five and
six-dimensional cases.

 \section{Applications}

 \subsection{The five-dimensional case}

This case is specific and can be also completely investigated:
$FLeib_5$ can be represented as a disjoint union of several
subsets:
$$FLeib_5=\emph{U}_{1}\bigcup \emph{U}_{2}\bigcup \emph{U}_{3}\bigcup
\emph{U}_{4}\bigcup\emph{U}_{5}\bigcup \emph{F},$$ where \\
$\emph{U}_{1}=\{L(\beta)\in FLeib_5:\beta _{3}\neq 0 \ \
\mbox{and} \ \ \gamma-2\beta_3^2\neq0 \},$ \\
$\emph{U}_{2}=\{L(\beta)\in FLeib_5:\beta _{3}\neq 0, \ \
\gamma-2\beta_3^2=0 \ \ \mbox{and} \ \ \beta_4\neq0 \},$ \\
$\emph{U}_{3}=\{L(\beta)\in FLeib_5:\beta _{3}\neq 0,  \ \
\gamma-2\beta_3^2=0 \ \ \mbox{and} \ \ \beta_4=0 \},$ \\
$\emph{U}_{4}=\{L(\beta)\in FLeib_5:\beta _{3}= 0, \ \ \gamma\neq0
\},$ \\ $\emph{U}_{5}=\{L(\beta)\in FLeib_5:\beta
_{3}= 0, \ \ \gamma=0 \ \ \mbox{and} \ \ \beta_4\neq0 \},$ \\
$\emph{F}=\{L(\beta)\in FLeib_5:\beta _{3}= 0,  \ \ \gamma=0 \ \
\mbox{and} \ \ \beta_4=0 \}.$

\textbf{Proposition 4.} $i)$ Two algebras $L(\beta)$ and
$L(\beta')$ from $\emph{U}_1$ are isomorphic if and only if
$$ \frac{\gamma}{\beta_{3}^{2}}=\frac{\gamma'}{\beta_{3}'^{2}}.$$

$ii)$ For any $\lambda \in \textbf{C}$ there is an algebra
$L(\beta)$ from $\emph{U}$ such that $
\frac{\gamma}{\beta_{3}^{2}}=\lambda.$

\textbf{Proposition 5.}

$\emph{a}$) The algebras from the set $U_2$ are isomorphic to the
algebra $L(1,1,2);$

$\emph{b}$) The algebras from the set $U_3$ are isomorphic to the
algebra $L(1,0,2);$

$\emph{c}$) The algebras from the set $U_4$ are isomorphic to the
algebra $L(0,0,1);$

$\emph{d}$) The algebras from the set $U_5$ are isomorphic to the
algebra $L(0,1,0);$

$\emph{e}$) The algebras from the set $F$ are isomorphic to the
algebra $L(0,0,0).$

\textbf{Theorem 6.}

Any 5-dimensional complex filiform Leibniz algebra from $FLeib_5$
is isomorphic to one of the following pairwise nonisomorphic
non-Lie filiform complex Leibniz algebras
$L=<e_1,e_2,e_3,e_4,e_5>$ whose commutation relations are (omitted
products are assumed to be zero):

1)$L(1,0,\lambda):$

\ \ $\left\{\begin{array}{ll} [e_1,e_1]=e_{3},   [e_3,e_1]=e_{4},
[e_4,e_1]=e_{5}, [e_1,e_2]=e_4, & \\[1mm]
[e_2,e_2]=\lambda e_5, [e_3,e_2]=e_{5}, \mbox{where} \ \ \lambda \in \textbf{C}.\\
\end{array}
\right.$

2)$L(1,1,2):$

\ \ $\left\{\begin{array}{ll} [e_1,e_1]=e_{3},   [e_3,e_1]=e_{4},
[e_4,e_1]=e_{5}, [e_1,e_2]=e_4+e_5, & \\[1mm]
[e_2,e_2]=2 e_5, [e_3,e_2]=e_{5}.\\
\end{array}
\right.$

3)$L(0,0,1):$

\ \ $\left\{\begin{array}{ll} [e_1,e_1]=e_{3},   [e_3,e_1]=e_{4},
[e_4,e_1]=e_{5},
[e_2,e_2]= e_5.& \\[1mm]
\end{array}
\right.$

4)$L(0,1,0):$

\ \ $\left\{\begin{array}{ll} [e_1,e_1]=e_{3},   [e_3,e_1]=e_{4},
[e_4,e_1]=e_{5}, [e_1,e_2]=e_5. & \\[1mm]
\end{array}
\right.$

5)$L(0,0,0):$

\ \ $\left\{\begin{array}{ll} [e_1,e_1]=e_{3},   [e_3,e_1]=e_{4},
[e_4,e_1]=e_{5} & \\[1mm]
\end{array}
\right.$

\subsection{The six-dimensional case}

The set $FLeib_6$ can be represented as a disjoint union of the
following subsets:
$$FLeib_6=\emph{U}_{1}\bigcup \emph{U}_{2}\bigcup \emph{U}_{3}\bigcup
\emph{U}_{4}\bigcup\emph{U}_{5}\bigcup
\emph{U}_{6}\bigcup\emph{U}_{7}\bigcup\emph{U}_{8}\bigcup\emph{U}_{9}\bigcup
\emph{U}_{10}\bigcup\emph{U}_{11}\bigcup\emph{F},$$ where \\
$\emph{U}_{1}=\{L(\beta)\in FLeib_6:\beta _{3}\neq 0, \ \
\beta_4\neq0 \ \ \mbox{and} \ \ \gamma \neq0 \},$ \\
$\emph{U}_{2}=\{L(\beta)\in FLeib_6:\beta _{3}\neq 0, \ \
\beta_4\neq0,  \ \ \gamma=0 \ \ \mbox{and} \ \
\beta_3\beta_4\neq4\beta_4^2 \},$ \\ $\emph{U}_{3}=\{L(\beta)\in
FLeib_6:\beta _{3}\neq 0 \ \ \beta_4\neq0,  \ \ \gamma=0 \ \
\mbox{and} \ \ \beta_3\beta_4=4\beta_4^2 \},$ \\
$\emph{U}_{4}=\{L(\beta)\in FLeib_6:\beta _{3}\neq 0, \ \
\beta_4=0,  \ \ \gamma\neq0 \},$ \\ $\emph{U}_{5}=\{L(\beta)\in
FLeib_6:\beta _{3}\neq 0, \ \ \beta_4=0,  \ \ \gamma=0 \},$
\\ $\emph{U}_{6}=\{L(\beta)\in FLeib_6:\beta _{3}= 0, \ \
\beta_4\neq0,  \ \ \gamma\neq0 \},$ \\ $\emph{U}_{7}=\{L(\beta)\in
FLeib_6:\beta _{3}= 0, \ \ \beta_4\neq0,  \ \ \gamma=0, \ \
\mbox{and} \ \ \beta_5\neq0 \},$ \\ $\emph{U}_{8}=\{L(\beta)\in
FLeib_6:\beta _{3}= 0, \ \ \beta_4\neq0,  \ \ \gamma=0, \ \
\mbox{and} \ \ \beta_5=0 \},$ \\ $\emph{U}_{9}=\{L(\beta)\in
FLeib_6:\beta _{3}= 0, \ \ \beta_4=0,  \ \ \beta_5\neq0, \ \
\mbox{and}
 \ \ \gamma\neq0 \},$ \\ $\emph{U}_{10}=\{L(\beta)\in
FLeib_6:\beta _{3}= 0, \ \ \beta_4=0,  \ \ \beta_5\neq0, \ \
\mbox{and}
 \ \ \gamma=0 \},$ \\ $\emph{U}_{11}=\{L(\beta)\in
FLeib_6:\beta _{3}= 0, \ \ \beta_4=0,  \ \ \beta_5=0, \ \
\mbox{and}
 \ \ \gamma\neq0 \},$ \\ $\emph{F}=\{L(\beta)\in
FLeib_6:\beta _{3}= 0, \ \ \beta_4=0,  \ \ \beta_5=0, \ \
\mbox{and}
 \ \ \gamma=0 \},$

\textbf{Proposition 7.} $i)$ Two algebras $L(\beta)$ and
$L(\beta')$ from $\emph{U}_1$ are isomorphic if and only if
$$ \varrho_3=\frac{2\beta_{3}\beta_{4}\gamma+4\beta_{3}^3\beta_{5}-5\beta_{3}^{2}\beta_{4}^2}{\gamma^2}=
\frac{2\beta_{3}'\beta_{4}'\gamma'+4\beta_{3}'^3\beta_{5}'-5\beta_{3}'^{2}\beta_{4}'^2}{\gamma'^2}.$$

$ii)$ For any $\lambda \in \textbf{C}$ there is an algebra
$L(\beta)$ from $\emph{U}_1$ such that
$\varrho_3=\frac{2\beta_{3}\beta_{4}\gamma+4\beta_{3}^3\beta_{5}-5\beta_{3}^{2}\beta_{4}^2}{\gamma^2}=\lambda.$

\textbf{Proposition 8.} $i)$ Two algebras $L(\beta)$ and
$L(\beta')$ from $\emph{U}_4$ are isomorphic if and only if
$$ \varrho_3=\frac{4\beta_{3}^3\beta_{5}}{\gamma^2}=\frac{4\beta_{3}'^3\beta_{5}'}{\gamma'^2}.$$

$ii)$ For any $\lambda \in \textbf{C}$ there is an algebra
$L(\beta)$ from $\emph{U}_4$ such that
$\varrho_3=\frac{4\beta_{3}^3\beta_{5}}{\gamma^2}=\lambda.$

\textbf{Proposition 9.}

$\emph{a}$) The algebras from the set $U_2$ are isomorphic to the
algebra $L(1,0,1,0);$

$\emph{b}$) The algebras from the set $U_3$ are isomorphic to the
algebra $L(1,0,0,0);$

$\emph{c}$) The algebras from the set $U_5$ are isomorphic to the
algebra $L(1,1,0,0);$

$\emph{d}$) The algebras from the set $U_6$ are isomorphic to the
algebra $L(0,1,0,1);$

$\emph{e}$) The algebras from the set $U_7$ are isomorphic to the
algebra $L(0,1,1,0);$

$\emph{f}$) The algebras from the set $U_8$ are isomorphic to the
algebra $L(0,1,0,0);$

$\emph{g}$) The algebras from the set $U_9$ are isomorphic to the
algebra $L(0,0,0,1);$

$\emph{h}$) The algebras from the set $U_{10}$ are isomorphic to
the algebra $L(0,0,1,0);$

$\emph{k}$) The algebras from the set $U_{11}$ are isomorphic to
the algebra $L(0,0,1,1);$

$\emph{l}$) The algebras from the set $F$ are isomorphic to the
algebra $L(0,0,0,0).$

\textbf{Theorem 10.}

Any 6-dimensional complex filiform Leibniz algebra from $FLeib_6$
is isomorphic to one of the following pairwise nonisomorphic
non-Lie filiform complex Leibniz algebras
$L=<e_1,e_2,e_3,e_4,e_5,e_6>$ whose commutation relations
are(omitted products are supposed to be zero):

1)$L(1,0,\lambda,1):$

\ \ $\left\{\begin{array}{ll} [e_1,e_1]=e_{3},   [e_3,e_1]=e_{4},
[e_4,e_1]=e_{5}, [e_5,e_1]=e_{6}, [e_1,e_2]=e_4+\lambda e_6, & \\[1mm]
[e_2,e_2]=e_6, [e_3,e_2]=e_{5}, [e_4,e_2]=e_{6},\mbox{where} \ \ \lambda \in \textbf{C}.\\
\end{array}
\right.$

2)$L(1,0,1,0):$

\ \ $\left\{\begin{array}{ll} [e_1,e_1]=e_{3},   [e_3,e_1]=e_{4},
[e_4,e_1]=e_{5}, [e_5,e_1]=e_{6}, [e_1,e_2]=e_4+e_6, & \\[1mm]
[e_3,e_2]=e_{5}, [e_4,e_2]=e_{6}.\\
\end{array}
\right.$

3)$L(1,0,0,0):$

\ \ $\left\{\begin{array}{ll} [e_1,e_1]=e_{3},   [e_3,e_1]=e_{4},
[e_4,e_1]=e_{5}, [e_5,e_1]=e_{6}, [e_1,e_2]=e_4, & \\[1mm]
[e_3,e_2]=e_{5}, [e_4,e_2]=e_{6}.\\
\end{array}
\right.$

4)$L(1,1,\lambda,1):$

\ \ $\left\{\begin{array}{ll} [e_1,e_1]=e_{3},   [e_3,e_1]=e_{4},
[e_4,e_1]=e_{5}, [e_5,e_1]=e_{6}, [e_1,e_2]=e_4+e_5+\lambda e_6, & \\[1mm]
[e_2,e_2]=e_{6}, [e_3,e_2]=e_{5}+e_6, [e_4,e_2]=e_{6}, \ \ \lambda \in \textbf{C}.\\
\end{array}
\right.$

5)$L(1,1,0,0):$

\ \ $\left\{\begin{array}{ll} [e_1,e_1]=e_{3},   [e_3,e_1]=e_{4},
[e_4,e_1]=e_{5}, [e_5,e_1]=e_{6}, [e_1,e_2]=e_4+e_5, & \\[1mm]
[e_3,e_2]=e_{5}+e_6, [e_4,e_2]=e_{6}.\\
\end{array}
\right.$

6)$L(0,1,0,1):$

\ \ $\left\{\begin{array}{ll} [e_1,e_1]=e_{3},   [e_3,e_1]=e_{4},
[e_4,e_1]=e_{5}, [e_5,e_1]=e_{6}, [e_1,e_2]=e_5, & \\[1mm]
[e_2,e_2]=e_6,[e_3,e_2]=e_6.\\
\end{array}
\right.$

7)$L(0,1,1,0):$

\ \ $\left\{\begin{array}{ll} [e_1,e_1]=e_{3},   [e_3,e_1]=e_{4},
[e_4,e_1]=e_{5}, [e_5,e_1]=e_{6}, [e_1,e_2]=e_5+e_6, & \\[1mm]
[e_3,e_2]=e_6.\\
\end{array}
\right.$

8)$L(0,1,0,0):$

\ \ $\left\{\begin{array}{ll} [e_1,e_1]=e_{3},   [e_3,e_1]=e_{4},
[e_4,e_1]=e_{5}, [e_5,e_1]=e_{6}, [e_1,e_2]=e_5, & \\[1mm]
[e_3,e_2]=e_6.\\
\end{array}
\right.$

9)$L(0,0,0,1):$

\ \ $\left\{\begin{array}{ll} [e_1,e_1]=e_{3},   [e_3,e_1]=e_{4},
[e_4,e_1]=e_{5}, [e_5,e_1]=e_{6}, [e_2,e_2]=e_6.
\\
\end{array}
\right.$

10)$L(0,0,1,0):$

\ \ $\left\{\begin{array}{ll} [e_1,e_1]=e_{3},   [e_3,e_1]=e_{4},
[e_4,e_1]=e_{5}, [e_5,e_1]=e_{6}, [e_1,e_2]=e_6.
\\
\end{array}
\right.$

11)$L(0,0,1,1):$

\ \ $\left\{\begin{array}{ll} [e_1,e_1]=e_{3},   [e_3,e_1]=e_{4},
[e_4,e_1]=e_{5}, [e_5,e_1]=e_{6}, [e_1,e_2]=e_6, [e_2,e_2]=e_6.
\\
\end{array}
\right.$

12)$L(0,0,0,0):$

\ \ $\left\{\begin{array}{ll} [e_1,e_1]=e_{3},   [e_3,e_1]=e_{4},
[e_4,e_1]=e_{5}, [e_5,e_1]=e_{6}.
\\
\end{array}
\right.$

\end{document}